\documentclass[12pt,a4paper]{amsart}
\usepackage[foot]{amsaddr}

\makeatletter
\renewcommand\normalsize{%
    \@setfontsize\normalsize{11.7}{14pt plus .3pt minus .3pt}%
    \abovedisplayskip 10\p@ \@plus4\p@ \@minus4\p@
    \abovedisplayshortskip 6\p@ \@plus2\p@
    \belowdisplayshortskip 6\p@ \@plus2\p@
    \belowdisplayskip \abovedisplayskip}
\renewcommand\small{%
    \@setfontsize\small{9.5}{12\p@ plus .2\p@ minus .2\p@}%
    \abovedisplayskip 8.5\p@ \@plus4\p@ \@minus1\p@
    \belowdisplayskip \abovedisplayskip
    \abovedisplayshortskip \abovedisplayskip
    \belowdisplayshortskip \abovedisplayskip}
\renewcommand\footnotesize{%
    \@setfontsize\footnotesize{8.5}{9.25\p@ plus .1pt minus .1pt}
    \abovedisplayskip 6\p@ \@plus4\p@ \@minus1\p@
    \belowdisplayskip \abovedisplayskip
    \abovedisplayshortskip \abovedisplayskip
    \belowdisplayshortskip \abovedisplayskip}
\ifdefined\pdfpagewidth
\else
\fi
\calclayout
\makeatother

\usepackage[backend=bibtex, 
sortcites=true,
firstinits=true,
hyperref,
maxbibnames=99,
]{biblatex}

\bibliography{Bibliography}{}

\usepackage{enumerate}

\usepackage[centertags]{amsmath}
\usepackage{amsfonts}
\usepackage{amssymb}
\usepackage{amsthm}
\usepackage{newlfont}
\usepackage{mathtools}
\usepackage{tikz}

\theoremstyle{plain}
\newtheorem{theorem}{Theorem}[section]

\newtheorem{conjecture}[theorem]{Conjecture}

\theoremstyle{definition}

\newcommand{\R}{\mathbb{R}}

\newcommand{\diam}{\operatorname{diam}}

\newcommand{\Lip}{\operatorname{Lip}}
\newcommand{\capa}{\operatorname{cap}}

\newcommand{\E}{\mathcal{E}}
\newcommand{\EE}{\mathbb{E}}

\newcommand{\amod}{\alpha_{\operatorname{mod}}}

\newcommand{\Hidden}[1]{}

\DeclareMathOperator{\Ent}{Ent}

\usepackage{booktabs}
\usepackage{hyperref}

\begin{document}

\title{A counterexample to a conjecture by Salez and Youssef}
\author{
Florentin M\"unch}
\address[Florentin M\"unch]{Leipzig University, 04109 Leipzig, Germany}
\email{cfmuench@gmail.com}

\date{\today}
\maketitle

\begin{abstract}
Remarkable progress has been made in recent years to establish log-Sobolev type inequalities under the assumption of discrete Ricci curvature bounds.
More specfically, Salez and Youssef have proven that the log-Sobolev constant can be lower bounded by the Bakry Emery curvature lower bound divided by the logarithm of the sparsity parameter.
They conjectured that the same holds true when replacing Bakry Emery by Ollivier curvature which is often times easier to compute in practice.
In this paper, we show that this conjecture is wrong by giving a counter example on birth death chains of increasing length.
\end{abstract}



\section{Introduction}

The log-Sobolev constant $\alpha_{LSI}$ on a reversible Markov chain $(X,p)$ is an important invariant closely related to hypercontractivity and mixing time estimates, see the famous paper by Diaconis and Saloff-Coste \cite{diaconis1996logarithmic}. 
It is given by
\[
\alpha_{LSI} := \inf_{f>0} \frac{\E(\sqrt f)}{\Ent(f)}
\]
where
\[
\Ent(f) = \EE(f\log f) - \EE (f) \log \EE (f)
\]
and $\E (f) = \E(f,f) = \EE (\Gamma (f,g))$ with
\[
\Gamma (f,g)(x) := \frac 1 2 \sum_y p(x,y)(f(y)-f(x))(g(y)-g(x)).
\]
Here, $\EE$ denotes the expectation with respect to the invariant probability measure $\pi$ of the Markov chain. 
We remark that in the notation of \cite{salez2025intrinsic}, we have $\alpha_{LSI} = 1/t_{LS}$.
The log Sobolev inequality turns out to be equivalent to certain hypercontractivity properties of the heat semigroup $P_t = e^{t\Delta}$ with $\Delta: \R^X \to \R^X$ given by
\[
\Delta f(x) = \sum_y p(x,y)(f(y)-f(x)),
\]
see \cite{diaconis1996logarithmic}.
Particularly interesting is the relation of the log-Sobolev constant to various Ricci curvature notions of Markov chains.

The research in discrete Ricci curvature is a relatively young subject which recently gained significant momentum.
Ollivier and, independently, Joulin introduced a concept of Ricci curvature of Markov chains based on a 1-Wasserstein contraction property in 2009 \cite{ollivier2007ricci,ollivier2009ricci,
joulin2009new}. This Wasserstein contraction property is equivalent to the Dobrushin criterion already developed in the 1950s, but not linked yet to Ricci curvature, see \cite{dobrushin1956central}.

Another discrete Ricci curvature notion is based on the Bakry-Emery calculus introduced in \cite{BakryEmery85} which generalizes the Bochner formula for Riemannian manifolds to metric measure spaces. This curvature notion has first been applied to discrete settings in \cite{schmuckenschlager1998curvature} and \cite{lin2010ricci}. 

In 2003, Forman introduced a Ricci curvature notion for cell complexes via the Bochner Weitzenböck formula, connecting the Hodge to the Bochner Laplacian. It was shown in \cite{jost2021characterizations} that, when appropriately extending a graph to a 2 dimensional cell complex, the Forman curvature is locally equivalent to the Ollivier curvature (meaning that the appropriate extension can depend on the edge whose curvature is considered).

Finally, there is the entropic Ricci curvature introduced by Erbar and Maas, and independently, Mielke \cite{erbar2012ricci, mielke2013geodesic} based on convexity of the entropy along 2-Wasserstein geodesics where a modified version of the 2-Wasserstein distance.
As in their setting, the heat flow is the gradient flow of the entropy, the assumption of a lower Ricci bound, i.e., an assumption on the convexity of the entropy, implies an exponential decay of the entropy under the heat flow which is equivalent to a modified log-Sobolev inequality \cite{erbar2012ricci}.
Moreover, there is a long list of modifications of the above mentioned curvature notions \cite{bauer2015li,
samson2022entropic,rapaport2023criteria,
weber2023li,dier2017discrete,munch2020spectrally}.

Let us now discuss the connection between Bakry Emery and Ollivier Ricci curvature, and log-Sobolev inequalities. 
Perhaps most famously, there was the Peres-Tetali conjecture asking whether a modified log-Sobolev constant $\amod$ (arising from replacing $\E(\sqrt{f})$ by $\E(f,\log f)$ in the definition of $\alpha_{LSI}$, see \cite{bobkov2006modified}) can be lower bounded by a constant times the lower Ollivier Ricci curvature bound.
In \cite{munch2023ollivier}, this conjecture was disproved via birth death chains on three vertices only.
In the same paper, it was shown that the modified log-Sobolev is lower bounded by the minimum Ollivier Ricci curvature bound, when additionally assuming non-negative sectional curvature, i.e., $W_\infty(p(x,\cdot),p(y,\cdot)) \leq d(x,y)$ for all vertices $x,y$ where $W_\infty$ is the infinity Wasserstein distance. This result has been improved and applied to various examples in \cite{caputo2025entropy}, and gave new lower bounds on the modified log-Sobolev constant for e.g.  zero range process in the mean field case.

Let us briefly discuss the relation between the log-Sobolev constant $\alpha$ and its modified version $\amod$.
It is well known that
\[
4\alpha_{LSI} \leq \amod \leq 2\lambda
\]
 where $\lambda$ is the spectral gap of the Laplacian. Moreover, it was recently shown  by Salez and Youssef in
 \cite{salez2025intrinsic}, and a slightly weaker version in \cite{salez2023upgrading} that
 \[
 \amod \leq  15\alpha_{LSI} \log d 
 \]
where $$d := \max\{1/p(x,y): p(x,y)>0\}$$ is the sparsity constant.

We now return to the connection between discrete Ricci curvature and log-Sobolev inequalities.
In \cite{salez2025intrinsic}, it was shown that the log-Sobolev constant can be estimated via Bakry Emery curvature bounds.

\begin{theorem}[{\cite[Theorem~3]{salez2025intrinsic}}] \label{thmBElogSob}
Let $(X,p)$ be a reversible Markov chain with Bakry-Emery curvature bounded from below by $K>0$, i.e.,
\[
\Gamma P_t f \leq e^{-2Kt} P_t \Gamma f
\]
where $P_t$ is the associated heat semigroup, then
\begin{align*}
\alpha_{LSI} \geq \frac{K}{33 \log d}.
\end{align*}
\end{theorem}
This is a remarkable result as this proves the Peres-Tetali conjecture  $\amod \gtrsim K$ under the assumption
that the upper bound in $ \amod \leq  15\alpha_{LSI} \log d $ is attained or almost attained.

It is very natural to ask whether a similar estimate holds true when replacing the Bakry-Emery curvature by the Ollivier curvature, without the strong assumption of non-negative sectional curvature. This is the subject of the conjecture by Salez and Youssef.
\begin{conjecture}[{\cite[Conjecture~1]{salez2025intrinsic}}]
\label{con:SalezYoussef}
Let $(X,p)$ be a reversible Markov chain with Ollivier curvature bounded from below by $K>0$, i.e.,
\[
\Lip(P_t f) \leq e^{-Kt}\Lip(f)\quad  \mbox{ for all } f \in \R^X \mbox{ and } t>0,
\]
then
\begin{align*}
\alpha_{LSI} \geq c\cdot \frac{K}{\log d}
\end{align*}
where $c$ is a universal constant, and $\Lip$ is the Lipschitz constant with respect to the combinatorial distance (also called hop count distance).
\end{conjecture}
Here, we formulated the Ollivier curvature lower bound by its semigroup characterization to demonstrate the conceptual similarity between Ollivier and Bakry Emery curvature.
Indeed, the Ollivier curvature $\kappa$ is a local quantity and is classically defined as
\[
\kappa(x,y) = 1-W_1(p(x,\cdot),p(y,\cdot))
\]
for neighbors $x,y$, where $W_1$ is the 1-Wasserstein distance given by 
\[
W_1(\mu,\nu) = \sup_{f} \int f d\mu - \int f d\nu 
\]
where the infimum is taken over all 1-Lipschitz functions.
It was shown in \cite{joulin2009new} and \cite{munch2017ollivier} that, assuming laziness, a lower bound of the Ollivier curvature, i.e., $\kappa(x,y) \geq K$ for all neighbors $x,y$, is equivalent to the Lipschitz contraction of the heat semigroup with rate $K$, i.e., $\Lip(P_t f) \leq e^{-Kt} \Lip(f)$.

In this note, we will give a counter example to the above conjecture by constructing suitable birth death chains.

\section{A counter example to the Salez-Youssef conjecture}

For estimating the log-Sobolev constant on birth death chains, its isocapacitary characterization is particularly useful as the isocapacitary profile is very easy to compute on birth death chains.
The following characterization was established in \cite{schlichting2019poincare} in the discrete setting.
\begin{theorem}\label{thm:LScapacity}
Let $(X,p)$ be a reversible Markov chain.
There exist global constants $c,C$ such that
\[
c\alpha_{LSI} \leq \inf_{A,B,\pi(A)\geq \frac 1 2} \frac{\capa(A,B)}{\pi(B)|\log \pi(B)|} \leq C\alpha_{LSI}
\]
where the capacity for $A,B \subseteq X$ is defined as
\[
\capa(A,B) := \inf \{\E(f):f|_A=0, f|_B=1\}
\]
and $\mathcal E$ is the Dirichlet energy.
\end{theorem}
We notice that the optimizer $f$ for the capacity always exists and satisfies $\Delta f = 0$ outside of $A$ and $B$ by a variational principle.

A particularly nice feature of birth death chains is that the capacity can be calculated very easily:
\[
\frac 1 {\capa(\{1,\ldots,a\}, \{b,\ldots, |X|\})} = 
\sum_{k=a}^{b-1} \frac{1}{\pi(k)p(k,k+1)}
\]
for $a < b$, as $\frac{1}{\capa(A,B)}$ is precisely the effective resistance between $A$ and $B$ where each edge $(x,y)$ has resistance $\frac{1}{\pi(x)p(x,y)}$, and on a birth death chain, the edge resistances form a serial circuit and thus simply add up to the effective resistance.

Moreover, the Ollivier curvature of birth death chains on $\{0,\ldots,n\}$ is easy to compute. That is, the Ollivier curvature is bounded below from $K$ if and only if
\[
K \leq p(k,k+1) - p(k,k-1) - p(k+1,k+2) + p(k+1,k)
\]
for all $k \in \{0,\ldots,n-1\}$ where by convention, $p(0,-1)=0=p(n,n+1)$, see e.g. \cite[Theorem~4.3]{joulin2009new}.

We now construct a sequence of birth death chains as a counter example to Conjecture~\ref{con:SalezYoussef}.
We first discuss a few obstructions.
\begin{itemize}
\item If $p(k,k+1)$ is decreasing and $p(k,k-1)$ is increasing, then the Markov chain has non-negative sectional curvature, giving a lower bound on $\amod$ by \cite{munch2023ollivier} which implies the conjecture.
\item If the invariant measure is log-concave, then a lower Ollivier curvature bound $K$ implies a lower Bakry Emery curvature bound $\frac 1 2 K$, see \cite[Theorem~3]{hua2023ricci}  and the conjecture follows from Theorem~\ref{thmBElogSob}.
\end{itemize}
Hence, for our counter example, we must avoid these cases.
We first give a meta-description of these birth death chains.
\begin{itemize}
\item They have $3n$ vertices and Ollivier Ricci curvature $\frac{1}{4n^2}$
\item On the first $n$ vertices, the jump rates are small (i.e. the forward jump rates of order $1/n^2$ and the backward jump rates of order $1/n$) and the invariant measure decays roughly exponentially in space.
\item On the latter $2n$ vertices, the jump rates are large (of constant order) and the invariant measure stays roughly constant.
\end{itemize}
Specifically, we define
\begin{align*}
4p(k,k+1) &= \begin{cases} \frac 1 {n^2} &: 1 \leq k \leq n \\
1 -\frac 1 n - \frac{k}{n^2}&: n< k \leq 3n-1,
\end{cases}\\
4p(k,k-1) &= \begin{cases} \frac 1 n + \frac {k+1} {n^2} &: 2 \leq k \leq n \\
1 &: n< k \leq 3n.
\end{cases}
\end{align*}
We notice that for $2 \leq k \leq 3n-1$, we have
\[
4c_k:=  4p(k,k-1) - 4p(k,k+1) = \frac 1 n + \frac k {n^2}.
\]
The Ollivier curvature $\kappa(k,k+1)$ therefore satisfies
\[
4\kappa(k,k+1) = 4c_{k+1} - 4c_k = \frac 1 {n^2}
\]
for $2 \leq k \leq 3n-2$.
Moreover, for $n \geq 4$,
\begin{align*}
4\kappa(1,2) &= 4p(1,2) + 4p(2,1) - 4p(2,3) = \frac 1 n + \frac{3}{n^2} \geq \frac 1 {n^2}, \qquad \mbox{and}
\\ 4\kappa(3n-1,3n) &= 4p(3n,3n-1) + 4p(3n-1,3n) - 4p(3n-1,3n-2) \\&= 1 - \frac 1 n - \frac{3n-1}{n^2}
= 1 - \frac 4 n + \frac 1 {n^2} \geq  \frac 1 {n^2}.
\end{align*}
In summary, we obtain $\kappa \geq \frac{1}{4n^2}$.

We now want to estimate the log-Sobolev constant 
via the isocapacitary characterization.
More precisely, we will choose $A=\{1\}$ and $B=\{2n,\ldots,3n\}$, and we will calculate
\[
 \frac{\capa(A,B)}{\pi(B)|\log \pi(B)|} 
\]
which will provide an upper bound to the log-Sobolev constant, as soon as we have $\pi(A) \geq \frac 1 2$.
To this end, we calculate the invariant measure $\pi$.
We have
\[
\frac{\pi(k)}{\pi(k+1)} = \begin{cases}
n+k+2&: 1\leq k \leq n-1 \\
n^2 &:k = n \\
\frac{1}{1-\frac 1 n - \frac{k}{n^2}}&: n+1 \leq k \leq 3n-1
\end{cases}
\]
In particular, we get $\pi(1) \geq \frac 1 2$ for $n$ large enough.
Moreover, for $n+1 \leq k \leq 3n$ and for $n$ large,
\[
1 \geq \frac{\pi(k)}{\pi(n+1)} \geq \left(1 - \frac 4{n} \right)^{2n} \geq e^{-10}.
\]
Hence,
\[
\pi(\{2n,\ldots,3n\}) \geq e^{-10} \pi(n+1) \cdot n.
\]
We next calculate the capacity
\begin{align*}
\frac 1 {\capa(\{1\},\{2n,\ldots,3n\})} &= \sum_{k=2}^{2n} \frac{1}{\pi(k)p(k,k-1)} \\
& \geq \sum_{k=n+1}^{2n} \frac{1}{\pi(k)p(k,k-1)} \\
& \geq \frac{n}{\pi(n+1)}
\end{align*}
and therefore,
\[
\frac  {\capa(\{1\},\{2n,\ldots,3n\})}{\pi(\{2n,\ldots,3n\})} \leq \frac{e^{10}}{n^2}.
\]
We finally estimate for $n \geq 4$
\begin{align*}
\pi(\{2n,\ldots,3n\}) &\leq (n+1)\pi(n+1)\\&=(n+1)
\pi(1) \cdot \prod_{k=1}^n\frac{\pi(k+1)}{\pi(k)}
\leq (n+1)\pi(1) n^{-n} \leq n^{-n/2}.
\end{align*}
Hence by Theorem~\ref{thm:LScapacity},
\begin{align*}
\alpha_{LSI} \leq
C\frac  {\capa(\{1\},\{2n,\ldots,3n\})}{\pi(\{2n,\ldots,3n\}) |\log \pi(\{2n,\ldots,3n\})|} &\leq C\frac{e^{10}}{\frac 1 2 n^3 \log n}  \\
&= \frac{\widetilde C}{ n^3 \log n}.
\end{align*}
We recall that the sparsity parameter $d$ is $4n^2 \leq n^3$ and the curvature $K$ is $\frac 1 {4n^2}$. Thus,
\[
\alpha_{LSI} \leq 12 \widetilde C \cdot \frac{K}{n \log d}. 
\]
Taking the limit $n\to \infty$ disproves  Conjecture~\ref{con:SalezYoussef}.

\section{Acknowledgments}
In a previous version of the article, there was a positive result included, namely
\[
\alpha_{LSI} \geq C\frac{K}{\diam \cdot \log d}
\]
where $K$ is the Ollivier Ricci bound, $\diam$ the diameter and $C$ a universal constant. In a private communication, Justin Salez explained that this bound trivially followed from the Lichnerowicz estimate $\lambda \geq K$ where $\lambda$ is the smallest positive eigenvalue of the Laplacian. The author wants to thank Justin Salez for pointing this.

\printbibliography


\end{document}